\documentclass[11pt]{article}
\usepackage[margin=1in]{geometry}
\usepackage{amsfonts,amsmath,amssymb}
\usepackage{algpseudocode}

\PassOptionsToPackage{obeyspaces}{url}
\usepackage[colorlinks=true,citecolor=blue,urlcolor=blue,linkcolor=blue,bookmarksopen=true]{hyperref}
\usepackage{amsthm}
\usepackage{tikz}
\usetikzlibrary{calc}

\usepackage{etoolbox}
\patchcmd{\thebibliography}{\leftmargin\labelwidth}{\leftmargin\labelwidth\addtolength\itemsep{-0.1\baselineskip}}{}{}

\author{Christopher Cox\thanks{Department of Mathematics, Iowa State University, Ames, IA, USA. \texttt{cocox@iastate.edu}. Supported in part through NSF RTG Grant DMS-1839918.}
\and Ryan R.\ Martin\thanks{Department of Mathematics, Iowa State University, Ames, IA, USA. \texttt{rymartin@iastate.edu}. Supported in part through Simons Collaboration Grants \#353292 and \#709641.}}

\title{The maximum number of $10$- and $12$-cycles in a planar graph}
\date{}

\usepackage[nameinlink,sort]{cleveref}

\newtheorem{theorem}{Theorem}
\newtheorem{lemma}[theorem]{Lemma}
\newtheorem{corollary}[theorem]{Corollary}

\crefname{conj}{conjecture}{conjectures}

\crefname{claim}{claim}{claims}

\crefname{prop}{proposition}{propositions}

\theoremstyle{definition}

\newtheorem{defn}[theorem]{Definition}
\crefname{defn}{definition}{definitions}

\newtheorem{remark}[theorem]{Remark}
\crefname{remark}{remark}{remarks}

\newtheorem{question}[theorem]{Question}
\crefname{question}{question}{questions}

\crefname{enumi}{part}{parts}

\makeatletter
\DeclareRobustCommand{\crefnosort}[1]{%
    \begingroup\@cref@sortfalse\cref{#1}\endgroup
}
\DeclareRobustCommand{\Crefnosort}[1]{%
    \begingroup\@cref@sortfalse\Cref{#1}\endgroup
}
\makeatother

\newcommand*{\eqdef}{\stackrel{\mbox{\normalfont\tiny{def}}}{=}}        
\newcommand*{\abs}[1]{\lvert #1\rvert}                

\renewcommand*{\epsilon}{\varepsilon}       
\newcommand*{\cyc}[2]{\cp_{#1}({#2})}

\newcommand*{\plan}{\mathcal{P}}
\newcommand*{\optb}{\beta}

\newcommand*{\sg}[1]{G_{#1}}

\DeclareMathOperator{\sgn}{sgn}
\DeclareMathOperator{\Opt}{Opt}
\DeclareMathOperator{\numb}{\mathbf{N}}
\DeclareMathOperator{\cp}{\mathbf{C}}
\newcommand*{\tsquare}{\mathbin{\square}}

\DeclareMathOperator{\supp}{supp}   

\begin{document}
\maketitle

\begin{abstract}
    For a fixed planar graph $H$, let $\numb_\plan(n,H)$ denote the maximum number of copies of $H$ in an $n$-vertex planar graph.
    In the case when $H$ is a cycle, the asymptotic value of $\numb_\plan(n,C_m)$ is currently known for $m\in\{3,4,5,6,8\}$.
    In this note, we extend this list by establishing $\numb_\plan(n,C_{10})\sim(n/5)^5$ and $\numb_\plan(n,C_{12})\sim(n/6)^6$.
    We prove this by answering the following question for $m\in\{5,6\}$, which is interesting in its own right: which probability mass $\mu$ on the edges of some clique maximizes the probability that $m$ independent samples from $\mu$ form an $m$-cycle?
\end{abstract}

\section{Introduction}

For graphs $G$ and $H$, let $\numb(G,H)$ denote the number of (unlabeled) copies of $H$ in $G$.
Furthermore, for a planar graph $H$, define
\[
    \numb_\plan(n,H)\eqdef\max\bigl\{\numb(G,H):G \text{ is an $n$-vertex planar graph}\bigr\}.
\]

The study of $\numb_\plan(n,H)$ was initiated by Hakimi and Schmeichel~\cite{hakimi_cycles}, who determined $\numb_\plan(n,C_3)$ and $\numb_\plan(n,C_4)$ exactly.
Alon and Caro~\cite{alon_biclique} continued this study by determining $\numb_\plan(n,K_{2,k})$ exactly for all $k$; in particular, they determined $\numb_\plan(n,P_3)$.
Gy\H{o}ri et al.~\cite{gyori_p4} later gave the exact value for $\numb_\plan(n,P_4)$, and the same authors determined $\numb_\plan(n,C_5)$ in~\cite{gyori_c5}.
Afterward, Ghosh et al.~\cite{ghosh_planarp5} asymptotically determined $\numb_\plan(n,P_5)$, and, very recently, the current authors~\cite{cox_paths} computed $\numb_\plan(n,P_7)$ asymptotically.

In~\cite{cox_paths} a general technique was introduced which allows one to bound $\numb_\plan(n,H)$ whenever $H$ exhibits a particular subdivision structure.
Using this technique, the authors established $\numb_\plan(n,C_6)\sim(n/3)^3$ and $\numb_\plan(n,C_8)\sim(n/4)^4$.
They furthermore conjectured that
\[
    \numb_\plan(n,C_{2m})=\biggl({n\over m}\biggr)^m+o(n^m)\qquad\text{for all }m\geq 3,
\]
and exhibited graphs meeting this lower bound.\footnote{It would not be surprising if this conjecture was suggested earlier, though it does not appear to be explicitly stated anywhere else in the literature. The closest reference we could find is~\cite{gyori_turan}, in which the lower-bound construction is mentioned.}
Currently, the best known general upper bound is
\[
    \numb_\plan(n,C_{2m})\leq {n^m\over m!}+o(n^m)\qquad\text{for }m\geq 5.
\]
\medskip

In this note, we make progress toward this conjecture by establishing the next two open cases:
\begin{theorem}\label{thm:c10c12}
    For $m\in\{5,6\}$,
    \[
        \numb_\plan(n,C_{2m})=\biggl({n\over m}\biggr)^m+o(n^m).
    \]
\end{theorem}

The paper is organized as follows.
In \Cref{sec:prelim}, we recall the technique developed in \cite{cox_paths} and extract the key ingredients necessary for the proof of \Cref{thm:c10c12}.
The proof of \Cref{thm:c10c12} is then the topic of \Cref{sec:proof}.

\section{Preliminaries}\label{sec:prelim}

In \cite{cox_paths} it was shown that one can upper bound $\numb_\plan(n,C_{2m})$ for $m\geq 3$ by answering the following question, which is interesting in its own right:
\begin{question}
    Which probability mass $\mu$ on the edges of some clique maximizes the probability that $m$ independent samples from $\mu$ form a copy of $C_m$?
\end{question}

To formalize this question, we lay out the following definition.
For a graph $G$ and an integer $m\geq 3$, let $\cyc mG$ denote the set of (unlabeled) copies of $C_m$ in $G$, so $\numb(G,C_m)=\abs{\cyc mG}$.
Additionally, for a finite set $X$, let $K_X$ denote the clique on vertex-set $X$.

\begin{defn}
    Let $X$ be a finite set and let $\mu$ be a probability mass on ${X\choose 2}$.
    \begin{enumerate}
        \item For a subgraph $G\subseteq K_X$, define
            \[
                \mu(G)\eqdef\prod_{e\in E(G)}\mu(e).
            \]
        \item For an integer $m\geq 3$, define
            \begin{align*}
                \optb(\mu;m) &\eqdef \sum_{C\in\cyc m {K_X}}\mu(C),\qquad\text{and}\\
                \optb(m) &\eqdef \sup\biggl\{\optb(\mu;m):\text{$\mu$ a probability mass on ${X\choose 2}$ for some finite set $X$}\biggr\}.
            \end{align*}
    \end{enumerate}
\end{defn}

The quantity $\optb(m)$ yields an upper bound on $\numb_\plan(n,C_{2m})$:
\begin{theorem}[{\cite[Lemma 2.5]{cox_paths}}]\label{thm:reduction}
    For $m\geq 3$,
    \[
        \numb_\plan(n,C_{2m})\leq \optb(m)\cdot n^m+O(n^{m-1/5}).
    \]
\end{theorem}
The argument here establishes that, as $n\to\infty$, the extremal graphs for $\numb_\plan(n,C_{2m})$ may be approximated by taking some fixed graph $G$ and ``blowing up'' the edges into independent sets of various sizes.
The probability mass $\mu$ is a compact way to represent the relative sizes of each of these independent sets, and $\optb(\mu;m)$ is a normalized count of the number of $C_{2m}$'s in the resulting blow-up.
Interestingly, planarity plays only a minor role in this argument.

In \cite{cox_paths}, the authors conjectured that $\optb(m)=m^{-m}$, which is the value achieved by the uniform distribution on $E(C_m)$, and they established the cases of $m=3$ and $m=4$.
In this paper we prove that $\optb(m)=m^{-m}$ for $m\in\{5,6\}$, which will then imply \Cref{thm:c10c12}.
\medskip

In order to accomplish this goal, we will need to understand the structure of the probability masses at play.
The following definition lays out two important aspects of such a probability mass.

\begin{defn}
    Fix a finite set $X$ and let $\mu$ be a probability mass on ${X\choose 2}$.
    \begin{enumerate}
        \item For $x\in X$, define
            \[
                \bar\mu(x)\eqdef\sum_{y\in X\setminus\{x\}}\mu(xy),
            \]
            which is the probability that an edge sampled from $\mu$ is incident to $x$.
            It can also be thought of as the weighted degree of the vertex $x$.
            Note that $\sum_{x\in X}\bar\mu(x)=2$ thanks to the handshaking lemma.
        \item The \emph{support graph} of $\mu$ is the graph $\sg\mu$, which has $E(\sg\mu)=\supp\mu$ and $V(\sg\mu)=\supp\bar\mu$.
            Since $\sg\mu$ records the edges of positive mass under $\mu$, observe that
            \[
                \optb(\mu;m)=\sum_{C\in\cyc{m}{\sg\mu}}\mu(C).
            \]
    \end{enumerate}
\end{defn}

In~\cite[Corollary 4.7]{cox_paths}, it was shown that $\optb(m)$ is achieved for every $m\geq 3$, so we introduce the following notation:
\begin{defn}
    For an integer $m\geq 3$, denote by $\Opt(m)$ the set of all probability masses $\mu$ satisfying $\optb(\mu;m)=\optb(m)$.
\end{defn}

We will require structural results about such optimal masses which were established in \cite{cox_paths}.
\begin{lemma}[{\cite[Lemma 4.5]{cox_paths}}]\label{lem:regularity}
    Fix $m\geq 3$ and $\mu\in\Opt(m)$.
    Then,
    \begin{alignat*}{2}
        m\cdot\optb(\mu;m)\cdot\mu(e) &=\sum_{\substack{C\in\cyc m{\sg\mu}:\\ E(C)\ni e}}\mu(C),\qquad &&\text{for every } e\in\supp\mu,\quad\text{and}\\
        m\cdot\optb(\mu;m)\cdot\bar\mu(x) &= \sum_{\substack{C\in\cyc m{\sg\mu}:\\ V(C)\ni x}}2\cdot\mu(C),\qquad &&\text{for every }x\in\supp\bar\mu.
    \end{alignat*}
\end{lemma}

\begin{lemma}[{\cite[Lemma 4.6]{cox_paths}}]\label{lem:vertbound}
    Fix $m\geq 3$ and $\mu\in\Opt(m)$.
    If $z\in(0,1)$ satisfies
    \[
        1-{m\over 2}z>(1-z)^m,
    \]
    then $\bar\mu(x)>z$ for all $x\in\supp\bar\mu$.
\end{lemma}

Explicitly, we will employ the following two straightforward corollaries.
\begin{corollary}\label{cor:upperbound}
    Fix $m\geq 3$ and $\mu\in\Opt(m)$.
    Then,
    \begin{itemize}
        \item For any $e\in\supp\mu$, we have $\mu(e)\leq 1/m$ with equality if and only if $e\in E(C)$ for every $C\in\cyc m{\sg\mu}$, and
        \item For any $x\in\supp\bar\mu$, we have $\bar\mu(x)\leq 2/m$ with equality if and only if $x\in V(C)$ for every $C\in\cyc m{\sg\mu}$.
    \end{itemize}
\end{corollary}
\begin{proof}
    Using \Cref{lem:regularity}, we observe that
    \[
        m\cdot\optb(\mu;m)\cdot\mu(e)=\sum_{\substack{C\in\cyc m{\sg\mu}:\\ E(C)\ni e}}\mu(C)\leq\sum_{C\in\cyc m{\sg\mu}}\mu(C)=\optb(\mu;m),
    \]
    and hence the first claim follows.
    The proof of the second claim is analogous.
\end{proof}
\begin{corollary}\label{cor:rightstruct}
    Fix $m\in\{5,6\}$.
    If $\mu\in\Opt(m)$, then $\abs{\supp\bar\mu}=m$ and $\bar\mu(x)=2/m$ for every $x\in\supp\bar\mu$.
\end{corollary}
\begin{proof}
    Set $z=2/(m+1)$.
    For $m\in\{5,6\}$, it can be checked that $1-{m\over 2}z>(1-z)^m$.
    Thus, \Cref{lem:vertbound} implies that $\bar\mu(x)>z=2/(m+1)$ for every $x\in\supp\bar\mu$.
    From here, we see that
    \[
        2=\sum_{x\in\supp\bar\mu}\bar\mu(x)>{2\over m+1}\cdot\abs{\supp\bar\mu}\qquad\implies\qquad \abs{\supp\bar\mu}<m+1.
    \]
    As such, we know that $\abs{\supp\bar\mu}=m$, since certainly $\abs{\supp\bar\mu}\geq m$.
    Therefore, every copy of $C_m$ in $\sg\mu$ must use every vertex of $\supp\bar\mu$; hence $\bar\mu(x)=2/m$ for every $x\in\supp\bar\mu$ thanks to \Cref{cor:upperbound}.
\end{proof}

\Cref{cor:rightstruct} is the key observation which enables our arguments in the next section.
If one could extend \Cref{cor:rightstruct} to all $m\geq 7$, then approaching the full conjecture that $\optb(m)=m^{-m}$ would likely be significantly more tractable.

\begin{remark}
    For $\mu\in\Opt(7)$, \Cref{lem:vertbound} implies that $\bar\mu>0.246$ and so $\abs{\supp\bar\mu}\leq 8$.

    For $\mu\in\Opt(m)$ for general $m$, using the bound $1-z\leq e^{-z}$ and \Cref{lem:vertbound} yields $\bar\mu>1.593/m$ and so $\abs{\supp\bar\mu}< 1.256\cdot m$.
\end{remark}

\section{Proof of \Cref{thm:c10c12}}\label{sec:proof}

To prove \Cref{thm:c10c12}, it suffices to establish that $\optb(m)=m^{-m}$ for $m\in\{5,6\}$, thanks to \Cref{thm:reduction}.
\medskip

We begin by establishing an inequality on the edge-masses in an optimal probability mass.

\begin{lemma}\label{lem:edgeineq}
    Fix $m\in\{5,6\}$ and $\mu\in\Opt(m)$.
    For any $e\in\supp\mu$,
    \[
        \biggl({2\over 2-m\cdot\mu(e)}\biggr)^4\biggl({m-4\over m-4+m\cdot\mu(e)}\biggr)^{m-4}\bigl(1-m\cdot\mu(e)\bigr)\leq 1.
    \]
\end{lemma}
\begin{proof}
    Fix any $e\in\supp\mu$ and define
    \[
        a\eqdef{2\over 2-m\cdot\mu(e)},\qquad\text{and}\qquad b\eqdef{m-4\over m-4+m\cdot\mu(e)}.
    \]
    Observe that $a,b\geq 0$ since $m\geq 5$ and $\mu(e)\in(0,1/m]$ (\Cref{cor:rightstruct}).
    We define a new probability mass $\mu'$ by
    \[
        \mu'(s)\eqdef\begin{cases}
            0 & \text{if }s=e,\\
            a\cdot\mu(s) & \text{if }\abs{s\cap e}=1,\\
            b\cdot\mu(s) & \text{otherwise}.
        \end{cases}
    \]
    Suppose that $e=xy$.
    \Cref{cor:rightstruct} tells us that $\bar\mu(x)=\bar\mu(y)=2/m$, so
    \begin{align*}
        \sum_{s\in\supp\mu'}\mu'(s) &=a\cdot\sum_{\abs{s\cap e}=1}\mu(s)+b\cdot\sum_{s\cap e=\varnothing}\mu(s)\\
                                    &=a\cdot\bigl(\bar\mu(x)+\bar\mu(y)-2\mu(e)\bigr)+b\cdot\bigl(1-\bar\mu(x)-\bar\mu(y)+\mu(e)\bigr)\\
                                    &=a\cdot\biggl({4\over m}-2\mu(e)\biggr)+b\cdot\biggl(1-{4\over m}+\mu(e)\biggr)\\
                                    &={4\over m}+{m-4\over m}=1;
    \end{align*}
    therefore, $\mu'$ is indeed a probability mass.

    Since $\abs{\supp\bar\mu}=\abs{\supp\bar\mu'}=m$ (\Cref{cor:rightstruct}) and $e\notin\supp\mu'$, we know that any $C\in\cyc m{\sg{\mu'}}$ has exactly $4$ edges incident to the pair $\{x,y\}$.
    Furthermore, $\optb(\mu;m)\geq\optb(\mu';m)$ since $\mu\in\Opt(m)$.
    Thus, by appealing to \Cref{lem:regularity}, we compute
    \begin{align*}
        \optb(\mu;m)\geq\optb(\mu';m) &=a^4\cdot b^{m-4}\cdot\sum_{C\in\cyc m{\sg{\mu'}}}\mu(C)\\
                                      &=a^4\cdot b^{m-4}\cdot\biggl(\optb(\mu;m)-\sum_{\substack{C\in\cyc m{\sg\mu}:\\ E(C)\ni e}}\mu(C)\biggr)\\
                                      &=a^4\cdot b^{m-4}\cdot\bigl(\optb(\mu;m)-m\cdot\optb(\mu;m)\cdot\mu(e)\bigr).
    \end{align*}
    Dividing both sides of this inequality by $\optb(\mu;m)$ proves the lemma.
\end{proof}

\Cref{lem:edgeineq} allows us to place lower bounds on $\mu(e)$ for $e\in\supp\mu$.
\begin{corollary}\label{cor:edgebound}
    Fix $m\in\{5,6\}$ and $\mu\in\Opt(m)$.
    If $z\in(0,1)$ satisfies
    \[
        \biggl({2\over 2-z}\biggr)^4\biggl({m-4\over m-4+z}\biggr)^{m-4}(1-z)>1,
    \]
    then $\mu(e)>z/m$ for all $e\in\supp\mu$.
\end{corollary}
\begin{proof}
    For $z\in[0,1)$ define the function
    \[
        f(z)\eqdef \biggl({2\over 2-z}\biggr)^4\biggl({m-4\over m-4+z}\biggr)^{m-4}(1-z).
    \]
    Observe that $f(z)>0$ for all $z\in[0,1)$, so we may compute
    \begin{align*}
        {f'(z)\over f(z)}={d\over dz}\log f(z)={4\over 2-z}-{m-4\over m-4+z}-{1\over 1-z}={(1-m)z^2+2z\over (2-z)(1-z)(m-4+z)}
    \end{align*}
    Again, $f(z)>0$ for all $z\in[0,1)$, so, since $m\geq 5$,
    \[
        \sgn f'(z)=\sgn\bigl((1-m)z^2+2z\bigr),\qquad\text{for all }z\in[0,1).
    \]
    In particular, we see that $f'(z)\geq 0$ for all $0\leq z\leq {2\over m-1}$ and $f'(z)\leq 0$ for all ${2\over m-1}\leq z<1$.
    Since $f(0)=1$ and $f(1)=0$, this implies that the curves $y=f(z)$ and $y=1$ intersect at $0$ and at some unique $z^*\in(0,1)$.
    Furthermore, $f(z)>1$ for all $z\in(0,z^*)$ and $f(z)<1$ for all $z\in(z^*,1)$.
    \medskip

    Now that we have a better understanding of the function $f$, the claim follows quickly.
    Suppose that $z\in(0,1)$ satisfies $f(z)>1$.
    If we were to have $0<\mu(e)\leq z/m$, then $0<m\cdot\mu(e)\leq z$.
    From above, this would then imply that $f(m\cdot\mu(e))>1$ as well, contradicting \Cref{lem:edgeineq}.
\end{proof}

\begin{remark}
\Cref{lem:edgeineq} (and hence \Cref{cor:edgebound}) follows solely from the observations laid out in \Cref{cor:rightstruct}.
Therefore, if \Cref{cor:rightstruct} can be extended to $m\geq 7$, then so can \Cref{lem:edgeineq} and \Cref{cor:edgebound}.
\end{remark}

With these observations in hand, we can determine $\optb(5)$ and $\optb(6)$.
Firstly, for a graph $G$, let $\mu_G$ denote the uniform distribution on $E(G)$.

\begin{theorem}
    $\optb(5)=5^{-5}$.
\end{theorem}
\begin{proof}
    Fix any $\mu\in\Opt(5)$.
    Thanks to \Cref{cor:rightstruct}, we know that $\abs{\supp\bar\mu}=5$ and that $\bar\mu(x)=2/5$ for all $x\in\supp\bar\mu$.
    Furthermore, letting $z=2/3$ in \Cref{cor:edgebound}, we see that $\mu(e)>2/15$ for all $e\in\supp\mu$.

    Now, certainly $\delta(\sg\mu)\geq 2$ because $\sg\mu$ has a spanning copy of $C_5$.
    Furthermore, for any $x\in\supp\bar\mu$, we have
    \[
        {2\over 5}=\bar\mu(x)>{2\over 15}\deg(x)\qquad\implies\qquad\deg(x)<3.
    \]
    We conclude that $\sg\mu$ is $2$-regular, and so we must have $\sg\mu=C_5$.
    Applying the arithmetic--geometric mean (AM--GM) inequality then yields
    \[
        \optb(\mu;5)=\mu(\sg\mu)\leq\biggl({1\over 5}\sum_{e\in\supp\mu}\mu(e)\biggr)^5={1\over 5^5}.
    \]
    with equality if and only if $\mu=\mu_{C_5}$.
\end{proof}

The proof that $\optb(C_6)=6^{-6}$ is similar, albeit slightly more involved.

\begin{theorem}
    $\optb(C_6)=6^{-6}$.
\end{theorem}
\begin{proof}
    We begin by observing that
    \[
        \optb(6)\geq\optb(\mu_{C_6};6)={1\over 6^6}.
    \]
    Now, fix any $\mu\in\Opt(6)$.
    Thanks to \Cref{cor:rightstruct}, we know that $\abs{\supp\bar\mu}=6$ and that $\bar\mu(x)=1/3$ for all $x\in\supp\bar\mu$.
    Furthermore, letting $z=6/11$ in \Cref{cor:edgebound}, we see that
    \begin{equation}\label{eqn:edgec6}
        \mu(e)>{1\over 11},\qquad\text{for all }e\in\supp\mu.
    \end{equation}

    Now, certainly $\delta(\sg\mu)\geq 2$ because $\sg\mu$ has a spanning copy of $C_6$.
    Furthermore, thanks to \cref{eqn:edgec6}, for any $x\in\supp\bar\mu$,
    \[
        {1\over 3}=\bar\mu(x)>{1\over 11}\deg(x)\qquad\implies\qquad\deg(x)<{11\over 3}<4.
    \]
    Therefore, each vertex of $\sg\mu$ must have either degree $2$ or degree $3$.
    In fact, we claim that $\sg\mu$ is either $2$- or $3$-regular.
    Indeed, if this were not the case, then, since $\sg\mu$ is certainly connected, there would be two adjacent vertices $x,y$ with $\deg(x)=2$ and $\deg(y)=3$.
    Now, since $\deg(x)=2$, and $\sg\mu$ has $6$ vertices, every copy of $C_6$ in $\sg\mu$ must use the edge $xy$ and so $\mu(xy)=1/6$ thanks to \Cref{cor:upperbound}.
    But then, one of the other two edges incident to $y$ must have mass at most
    \[
        {\bar\mu(y)-1/6\over 2}={1/3-1/6\over 2}={1\over 12}<{1\over 11},
    \]
    contradicting \cref{eqn:edgec6}.
    \medskip

    Next, we claim that $\sg\mu$ must actually be $2$-regular.
    To prove this, suppose for the sake of contradiction that $\sg\mu$ is $3$-regular.

    To begin, since $\sg\mu$ is $3$-regular, and hence has $9$ edges, we may apply \cref{eqn:edgec6} to see that for any $C\in\cyc{6}{\sg\mu}$,
    \begin{equation}\label{eqn:cycbound}
        \sum_{e\in E(C)}\mu(e)=1-\sum_{e\in E(\sg\mu)\setminus E(C)}\mu(e)<1-{3\over 11}={8\over 11}.
    \end{equation}

    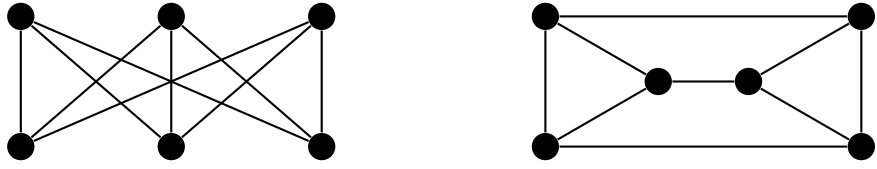
\begin{figure}[ht]
        \begin{center}
            \begin{tikzpicture}
                \foreach \x in {1,2,3} {
                    \node[circle, fill=black] (1\x) at (2*\x,-{sin(60)}) {};
                    \node[circle, fill=black] (2\x) at (2*\x,{sin(60)}) {};
                }
                \foreach \x in {1,2,3} {
                    \foreach \y in {1,2,3} {
                        \draw[thick] (1\x)--(2\y);
                    }
                }
            \end{tikzpicture}
            \hfil
            \begin{tikzpicture}
                \foreach \x in {1,2,3} {
                    \node[circle, fill=black] (1\x) at (120*\x:1) {};
                    \node[circle, fill=black] (2\x) at ($(3.2,0) - (-120*\x:1)$) {};
                    \draw[thick] (1\x)--(2\x);
                }
                \foreach \x in {1,2} {
                    \draw[thick] (\x1)--(\x2)--(\x3)--(\x1);
                }
            \end{tikzpicture}
        \end{center}
        \caption{The only two $3$-regular graphs on $6$ vertices: $K_{3,3}$ (left) and $K_3\tsquare K_2$ (right).\label{fig:k3k2}}
    \end{figure}
    Now, it is a routine exercise to show that the only $3$-regular graphs on $6$ vertices are $K_{3,3}$ and $K_3\square K_2$ (here, ``$\tsquare$'' denotes the Cartesian product of graphs; see \Cref{fig:k3k2} for a drawing of $K_3\tsquare K_2$).
    In either case, we have $\numb(G_\mu,C_6)\leq 6$; thus by applying the AM--GM inequality and \cref{eqn:cycbound}, we bound
    \begin{align*}
        {1\over 6^6}\leq \optb(6) = \optb(\mu;6) &= \sum_{C\in\cyc6{G_\mu}}\mu(C)\leq\sum_{C\in\cyc6{G_\mu}}\biggl({1\over 6}\sum_{e\in E(C)}\mu(e)\biggr)^6\\
                                                 &< 6\cdot\biggl({8\over 11}\biggr)^6\cdot{1\over 6^6}\leq {0.89\over 6^6}<{1\over 6^6};
    \end{align*}
    a contradiction.
    \medskip

    Therefore, we know that $\sg\mu$ is $2$-regular, and so $\sg\mu=C_6$.
    Applying the AM--GM inequality one final time then yields
    \[
        \optb(\mu;6)=\mu(\sg\mu)\leq\biggl({1\over 6}\sum_{e\in\supp\mu}\mu(e)\biggr)^6={1\over 6^6},
    \]
    with equality if and only if $\mu=\mu_{C_6}$.
    This concludes the proof.
\end{proof}

\bibliographystyle{abbrv}
\bibliography{references}

\end{document}